\documentclass[12pt]{article}
\usepackage{amsxtra,amssymb,amsthm,amsmath,latexsym}

\theoremstyle{plain}

\def\R{{\mathbb R}}

\def\oH{\buildrel\circ\over H}
\def\oH1{\buildrel\circ\over H\kern-.02in{}^1}

\def\d{\delta}

\def\be{\begin{equation}}
\def\ee{\end{equation}}

\begin{document}

%A.G.Ramm,  
%\begin{titlepage}

\title{ DSM for solving ill-conditioned linear algebraic systems 
   \thanks{key words: DSM (Dynamical Systems Method), ill-posed problems, 
ill-conditioned matrices
    }
}

\author{
A.G. Ramm\\
Mathematics Department, 
Kansas State University, \\
 Manhattan, KS 66506-2602, USA\\
ramm@math.ksu.edu\\
}

\date{}

\maketitle\thispagestyle{empty}

\begin{abstract} 
A standard way to solve linear algebraic systems $Au=f,\,\,(*)$ with 
ill-conditioned 
matrices $A$ is to use variational regularization. This leads to solving
the equation $(A^*A+aI)u=A^*f_\d$, where $a$ is a regularization 
parameter, and $f_\d$ are noisy data, $||f-f_\d||\leq \d$. 
Numerically it requires to calculate products of matrices $A^*A$ and 
inversion of the matrix $A^*A+aI$ which is also ill-conditioned if $a>0$ is 
small. We propose a new method for solving (*) stably, given noisy data 
$f_\d$. This method, the DSM (Dynamical Systems Method) is developed in 
this paper for selfadjoint $A$. It consists in solving a Cauchy problem
for systems of ordinary differential equations.

\end{abstract}

%\end{titlepage}

\section{Introduction} Consider a linear algebraic system
\be\label{e1.1} A^\ast Au=f, \quad f\in R(A), \ee where $A$
is a linear operator in $n-$dimensional Euclidean space,
$A=A^*$, $R(A)$ is the range of $A$, and $N:=\{u: Au=0\}$ is
the null-space of $A$. Let $k(A):=||A||||A^{-1}||$ denote
the condition number of $A$. If $A$ is singular, i.e., $N$
is not trivial, then we set $k(A)=\infty$. Problem (1.1) is
called ill-conditioned if $k(A)>>1$. In this case small
perturbations of $f$ may lead to large perturbations of the
solution $u$, so problem (1.1) is ill-posed. Such problems
are often solved by variational regularization. This method
consists in finding global minimizer of the functional
$F(u)=||Au-f_\d||^2+a||u||^2$, where $a=const>0$ is a
regularization parameter, and $f_\d$ are noisy data,
$||f-f_\d||\leq \d$. The global minimizer of the quadratic
functional $F$ is the unique solution to the linear
algebraic system $(A^*Au+aI)u_{a,\d}=A^*f_\d$, where $I$ is
the unit matrix. This system has a unique solution
$u_{a,\d}= (A^*Au+aI)^{-1}A^*f_\d$.  
Calculation of the matrix $A^*A$ requires multiplication of 
two matrices. For large $n$ this is a time-consuming 
operation. Condition number of the matrix $A^*A$ is 
$k^2(A)$, so it is much larger than $k(A)$ since $k(A)>>1.$
If $a$ is small, then
the condition number of the matrix $A^*A+aI$ is also large.
Therefore, inversion of the
matrix $A^*A+aI$ is numerically difficult if $a$ is small.
 An additional difficulty consists in
choosing the regularization parameter $a$
as a function of $\d$ in such a way that the element
$u_\d:=u_{a(\d),\d}$ would converge to a solution of (1.1)
as $\d\to 0$. There are a priori and a posteriori methods
for choosing such $a(\d)$. The theory of variational
regularization is presented in many books and papers (see
e.g., [1], Chapter 2).

Our goal is to develop a new method for for solving ill-conditioned 
problems (1.1). This method, which we call the DSM 
(Dynamical Systems Method),
does not require inversion of matrices and their multiplications.
It requires solving a Cauchy problem for a system of ordinary diferential
equations. The theoretical development of DSM is presented in [2],[3].
The author hopes that DSM will be an efficient numerical method for 
solving ill-conditioned linear algebraic systems.
The assumption that $A$ is selfadjoint, which is used for simplicity in 
this paper, can be relaxed. However, all generalizations will be 
considered elsewhere. Here we concentrate on the 
assumptions leading to the simplest arguments.

The idea of our method is simple: it is based on the formula
\be\label{e1.2} B^{-1}(e^{Bt}-I)=\int_0^te^{Bs}ds, \ee where
$B$ is a linear boundedly invertible operator. If
\be\label{e1.3} \lim_{t\to \infty}||e^{Bt}||=0, \ee then
\be\label{e1.4} B^{-1}=-\lim_{t\to \infty}\int_0^te^{Bs}ds.  
\ee The integral $\int_0^te^{B(t-s)}dsf$ is the solution to
the Cauchy problem \be\label{e1.5} \dot u=Bu+f,\quad u(0)=0,
\quad \dot u=\frac {du}{dt}. \ee Therefore, one can
calculate the inverse of an operator satisfying condition
(3) by solving a Cauchy problem.  We want to calculate
$A^{-1}f$. Let us take $B=i(A+iaI)$, where $I$ is the
identity operator and $a>0$ is a parameter which we take to
zero later. Condition (3) is satisfied if $A=A^*$ and
$a>0$ because under these assumptions $||e^{iAt}||=1$ and
$\lim_{t\to \infty}e^{-at}=0$. We write $A+ia$ in place of
$A+iaI$ below.
 
Consider the problem:
\be\label{e1.6} 
\dot u_a=i(A+ia)u_a+f, \quad u_a(0)=0. 
\ee
Its unique solution is $u_a=\int_0^t e^{i(A+ia)(t-s)}dsf$.
Our results are the following two theorems. 
 
{\bf Theorem 1.} {\it One has
\be\label{e1.7}
-i\lim_{a\to 0}\lim_{t\to \infty}u_a(t)=y,
\ee
where $y$ is the minimal-norm solution to (1.1).}

Note that if
$$\dot v_a=i(A+ia)v_a-if, \quad v_a(0)=0,
$$
then $\lim_{a\to 0}\lim_{t\to \infty}v_a(t)=y$.

If $f_\d$ is given in place of $f$, $||f-f_\d||\leq \d$,
then one solves the problem:
\be\label{e1.8}
\dot u_{a,\d}=i(A+ia)u_{a,\d}+f_\d, \quad u_{a,\d}(0)=0.
\ee
{\bf Theorem 2.} {\it If $t_\d$ and $a=a(\d)$ are such that
\be\label{e1.9}
\lim_{\d\to 0}t_\d=\infty,\quad \lim_{\d\to 0}a(\d)=0,\quad
\lim_{\d \to 0}\frac {\d}{a(\d)}=0,\quad \lim_{\d \to 
0}a(\d)t_\d=\infty,
\ee
then
\be\label{e1.10}
\lim_{\d\to 0}||u_\d-iy||=0,
\ee
where $u_\d:= u_{a(\d),\d}$.}

In the next Section proofs are given.

\section{Proofs}

{\bf Proof of Theorem 1.}  
By the argument in Section 1, we have
\be\label{e2.1}
u_{a}(t)=[-i(A+ia)]^{-1}(I-e^{i(A+ia)t})f,
\ee
and $\lim_{t\to \infty}||e^{i(A+ia)t}||=0$. Thus,
$\lim_{t\to \infty}u_{a}(t)=i(A+ia)^{-1}f$. Consequently,
$$-i\lim_{a\to 0}i(A+ia)^{-1}f=\lim_{a\to 
0}(A+ia)^{-1}Ay=y,$$
as claimed in Theorem 1.
Let us explain the last step. Using the spectral theorem for the 
selfadjoint operator $A$, one gets
\be\label{e2.2}
\lim_{a\to 0}||(A+ia)^{-1}Ay-y||^2=\lim_{a\to 
0}\int_{\R}|s(s+ia)^{-1}-1|^2d(E_sy,y)=\lim_{a\to 0}\int_{\R}\frac {a^2}
{a^2+s^2}d(E_sy,y)=0.
\ee
Here we have used the assumption $y\bot N$, which implies that
$\lim_{0<b\to 0}\int_{-b}^0d(E_sy,y)=0$. This relation 
allows one to pass
to the limit $a\to 0$ under the integral sign in (12).
Theorem 1 is proved. \hfill $\Box$

{\bf Remark 1.} In our case the operator $A$ is bounded,
so the integration in (12) is taken over a finite interval 
$[-||A||,||A||]$. Moreover, in a finite-dimensional space 
the spectrum of the operator $A$ consists of finitely many 
eigenvalues, and the integral in (12) reduces to a finite 
sum. Our proof of Theorem 1 allows one to use it
in infinite-dimensional spaces.  

{\bf Proof of Theorem 2.} One has
\be\label{e2.3}
||u_\d-y||\leq ||u_{a(\d)}-y||+||\int_0^{t_\d}e^{i(A+ia)(t-s)}ds(f_\d-f)||.
\ee
The argument given in the proof of Theorem 1 shows that
\be\label{e2.4}
\lim_{\d \to 0}||u_{a(\d)}(t_\d)-y||=0
\ee
provided that 
$$\lim_{\d \to 0}a(\d)=0,\,\, \lim_{\d \to 0}t_\d=\infty, 
\hbox {
and } \lim_{\d \to 0}a(\d)t_\d=\infty.$$
One estimates the integral
\be\label{e2.5}
J:=||\int_0^{t_\d}e^{i(A+ia)(t-s)}ds(f_\d-f)||\leq \frac {\d}{a}.
\ee
Thus
\be\label{e2.6}
||u_\d-iy||\leq  \frac {\d}{a} +o(1),
\ee
where the term $o(1)$ comes from (14).
Therefore, if $\lim_{\d \to 0} \frac {\d}{a(\d)}=0$ 
then $\lim_{\d \to 
0}||u_\d-y||=0$,
and Theorem 2 is proved. \hfill $\Box$

{\bf Remark 2.} The assumption $A=A^*$ allows one to give a 
short and 
simple proof of Theorems 1,2. However, this assumption can be replaced by 
more general assumptions. For example, one can assume that $A$ has 
Jordan chains of length one only. In other words, that the resolvent of 
$A$ has only simple poles. Under this more general assumption
the Cauchy problem we have used should be also modified, in general.
If $A=A^*$, then our proofs of theorems 1 and 2 remain valid for
operator $A$ in a Hilbert space, and not only in a finite-dimensional 
space.

\end{document}